%&amstex
% This is an AMS-TeX file and should be compiled
% using AMS-TeX.
% the command should be:  amstex file  or something
% like that, depending on what has been set up
% on the system you're using. 
\magnification=1200
\loadmsam
\loadmsbm
\loadeufm
\loadeusm
\UseAMSsymbols     
\input amssym.def
%\hsize=6.00 true in
%\hoffset=.50 true in
%\voffset=-0.1 true in
%\vsize=8.75 true in

\font\BIGtitle=cmr10 scaled\magstep3
\font\bigtitle=cmr10 scaled\magstep1
\font\boldsectionfont=cmb10 scaled\magstep1
\font\section=cmsy10 scaled\magstep1

\def\scr#1{{\fam\eusmfam\relax#1}}
\def\scrA{{\scr A}}
\def\scrB{{\scr B}}

\def\scrG{{\scr G}}
\def\scrH{{\scr H}}

\def\scrK{{\scr K}}
\def\scrJ{{\scr J}}

\def\scrM{{\scr M}}
\def\scrN{{\scr N}}
\def\scrO{{\scr O}}

\def\scrS{{\scr S}}

\def\scrR{{\scr R}}
\def\scrT{{\scr T}}

\def\db#1{{\fam\msbfam\relax#1}}

\def\dbA{{\db A}} 
\def\dbC{{\db C}} \def\dbD{{\db D}}
 \def\dbF{{\db F}}
\def\dbG{{\db G}} \def\dbH{{\db H}}

 \def\dbN{{\db N}}
 
\def\dbQ{{\db Q}} \def\dbR{{\db R}}

 \def\dbZ{{\db Z}}

\def\Rtil{\widetilde{R}}

\def\Ker{\text{Ker}}
\def\der{\text{der}}
\def\Sh{\hbox{\rm Sh}}

\def\ad{\text{ad}}

\def\Hom{\text{Hom}}
\def\End{\text{End}}

\def\Spec{\text{Spec}}

\def\Lie{\text{Lie}}

\def\leaderfill{\leaders\hbox to 1em
     {\hss.\hss}\hfill}
\def\nspace{\lineskip=1pt\baselineskip=12pt\lineskiplimit=0pt}

     %the way to use this is "\Proclaim{Theorem 1.1.}" for instance.

\def\finishproclaim{\par\rm
     \ifdim\lastskip<\medskipamount\removelastskip
     \penalty55\medskip\fi}

\def\references#1{\par
  \centerline{\boldsectionfont References}\smallskip
     \parindent=#1pt\nspace}
\def\Ref[#1]{\par\hang\indent\llap{\hbox to\parindent
     {[#1]\hfil\enspace}}\ignorespaces}
\def\Item#1{\par\smallskip\hang\indent\llap{\hbox to\parindent
     {#1\hfill$\,\,$}}\ignorespaces}
\def\ItemItem#1{\par\indent\hangindent2\parindent
     \hbox to \parindent{#1\hfill\enspace}\ignorespaces}

\def\Ge{{\mathchoice{\,{\scriptstyle\ge}\,}
  {\,{\scriptstyle\ge}\,}
  {\,{\scriptscriptstyle\ge}\,}{\,{\scriptscriptstyle\ge}\,}}}

\def\arrowsim{\,\smash{\mathop{\to}\limits^{\lower1.5pt
  \hbox{$\scriptstyle\sim$}}}\,}

\def\doublemaprights#1#2#3#4{\raise3pt\hbox{$\mathop{\,\,\hbox to
     #1pt{\rightarrowfill}\kern-30pt\lower3.95pt\hbox to
     #2pt{\rightarrowfill}\,\,}\limits_{#3}^{#4}$}}

\def\rightcapdownarrow{\raise9pt\hbox{$\ssize\cap$}\kern-7.75pt
     \Big\downarrow}

\def\rcapmapdown#1{\rightcapdownarrow\kern-1.0pt\vcenter{
     \hbox{$\scriptstyle#1$}}}

\def\rmapdown#1{\Big\downarrow\kern-1.0pt\vcenter{
     \hbox{$\scriptstyle#1$}}}
\def\rightsubsetarrow#1{{\ssize\subset}\kern-4.5pt\lower2.85pt
     \hbox to #1pt{\rightarrowfill}}
\def\longtwoheadedrightarrow#1{\raise2.2pt\hbox to #1pt{\hrulefill}
     \!\!\!\twoheadrightarrow}

\def\Hom{\operatorname{\hbox{Hom}}}

\NoBlackBoxes
\parindent=25pt
\document
\footline={\hfil}

\null
\vskip 0.4in
\centerline{\BIGtitle Three methods to prove the existence of integral}
\vskip 0.4in
\centerline{\BIGtitle canonical models of Shimura varieties of Hodge type} 
\vskip 0.2in 
\centerline{\bigtitle Adrian Vasiu, Binghamton University}
\vskip 0.2in
\centerline{November 18, 2008}
\footline={\hfill}
\null

\smallskip\noindent
{\bf ABSTRACT.} This is a survey of the three main methods developed in the last 15 years to prove the existence of integral canonical models of Shimura varieties of Hodge type. The only new part is formed by corrections to results of Kisin. 

\footline={\hss\tenrm \folio\hss}
\pageno=1

\bigskip\smallskip
\noindent
{\boldsectionfont 1. Introduction}

\medskip
Let $(G,X)$ be a Shimura variety of Hodge type. We fix an embedding $f:(G,X)\hookrightarrow (\pmb{GSp}(W,\psi),S)$ into a Shimura pair $(\pmb{GSp}(W,\psi),S)$ that defines a Siegel modular variety. Let $L$ be a $\dbZ$-lattice of $W$ such that we have a perfect alternating form $\psi:L\times L\to \dbZ$. Let $d\in\dbN$ be such that the rank of $L$ is $2d$. Let $\Sh(G,X)$ be the canonical model of $(G,X)$ over the reflex field $E(G,X)$.  

Let $p$ be a prime such that the reductive group $G_{\dbQ_p}$ extends to a reductive group scheme $G_{\dbZ_p}$ over $\dbZ_p$. It is known that there exists a reductive group scheme $G_{\dbZ_{(p)}}$ over $\dbZ_{(p)}$ whose generic fibre is $G$ and whose pull back to $\dbZ_p$ is $G_{\dbZ_p}$, cf. [Va1, Lemma 3.1.3]. We can choose $f$ such that the schematic closure of $G$ in $\pmb{GL}_{L\otimes_{\dbZ} \dbZ_{(p)}}$ is $G_{\dbZ_{(p)}}$, cf. [Va8, Part I, Lemma 4.2.1]. Let $v$ be a prime of the reflex field $E(G,X)$ that divides $p$; it is unramified over $p$. Let $O$ be the localization at $v$ of the ring of integers of $E(G,X)$. Let $k(v)$ be the residue field of $v$ (i.e., of $O$).  

Let $K_p:=\pmb{GSp}(L,\psi)(\dbZ_p)$ and let $H_p:=K_p\cap G(\dbQ_p)=G_{\dbZ_{(p)}}(\dbZ_p)$; they are hyperspecial subgroups of $\pmb{GSp}(L,\psi)(\dbQ_p)$ and $G(\dbQ_p)$ (respectively). The functorial morphism $\Sh(G,X)/H_p\to \Sh(\pmb{GSp}(W,\psi),S)_{E(G,X)}/K_p$ is a closed embedding, cf. [De1] and [De2]. Let $\scrM:=\text{proj.lim.}_{\dag\in\dbN\setminus p\dbN} \scrA_{d,1,\dag}$, where $\scrA_{d,1,\dag}$ is the Mumford's moduli scheme over $\dbZ_{(p)}$ that parametrizes isomorphism classes of principally polarized abelian schemes that have level-$\dag$ symplectic similitude structures. It is well known that one can identify $\scrM_{\dbQ}=\Sh(\pmb{GSp}(W,\psi),S)/K_p$ (for instance, cf. [De1]). Thus we can speak about the schematic closure $\scrN^{\text{cl}}$ of $\Sh(G,X)/H_p$ in $\scrM_O$ and about the normalization $\scrN$ of $\scrN^{\text{cl}}$. The goal of this survey is to report on the three main methods developed in the last 15 years to prove the following Theorem.

\bigskip\noindent
{\bf 1.1. Theorem.} {\it The $O$-scheme $\scrN$ is regular and formally smooth.}

\medskip
The three methods are presented in Sections 4 to 6. In Section 2 we include some complements on reductive group schemes. The basic notations to be used in order to detail on the three methods, are presented in Section 3.

\bigskip\noindent
{\bf 1.2. Remark.} The process of taking the normalization $\scrN$ of $\scrN^{\text{cl}}$ appears for the first time in [Va0] and thus in [Va1] and all its preliminary versions. In [Mi, Rms. 2.6 and 2.15] one uses only the schematic closure. The process of taking the normalization has been inserted due to weakness in the deformation theories available and due to the lack of progress in connection to the Tate and Hodge conjectures for abelian varieties. As pointed out in [Va1, Rm. 5.6.4], one does expect that always $\scrN=\scrN^{\text{cl}}$.  Our recent work on Langlands--Rapoport conjecture checks that in fact one has $\scrN=\scrN^{\text{cl}}$ in most cases.

\bigskip\smallskip
\noindent
{\boldsectionfont 2. Reductive group schemes}

\medskip
If $M$ is a free module of finite rank over a  commutative $\dbZ$-algebra $R$, then let $M^*:=\Hom(M,R)$, let $\scrT(M):=\oplus_{s,t\in\dbN\cup\{0\}} M^{\otimes s}\otimes_R M^{*\otimes t}$, and let $\scrS(M)$ be the $R$-module obtained from $M$ via the operations of direct sums, tensor products, duals, exterior powers, and symmetric powers. The definition of $\scrS(M)$ is too vague and for practical reasons one can take $\scrS(M):=\oplus_{n=1}^{\infty} \scrT(\Lambda^{\bullet}(S^{\bullet}(\scrT(M))))$, where $S^{\bullet}$ defines the symmetric poer algebra and $\Lambda^{\bullet}$ denotes the exterior power algebra.

It is known that there exists a family of tensors $(v_{\alpha})_{\alpha\in\scrJ}$ of $\scrT(W)$ such that $G$ is the subgroup of $\pmb{GL}_W$ that fixes $v_{\alpha}$ for all $\alpha\in\scrJ$ (for instance, see [De3, Prop. 3.1 c)]). One could ask if there exists an analogue result over $\dbZ_{(p)}$ instead of over $\dbQ$. Partial results were provided in the past as follows.

\bigskip
\noindent
{\bf 2.1. PEL type.} We assume that $f:(G,X)\hookrightarrow (\pmb{GSp}(W,\psi),S)$ is a PEL type embedding (i.e., $G$ is the identity component of the centralizer in $\pmb{GSp}(W,\psi)$ of a semisimple $\dbQ$--subalgebra of $\End(W)$). If $(G,X)$ is of either $A$ or $C$ type, then $G_{\dbZ_{(p)}}$ is the subgroup scheme of $\pmb{GL}_{L\otimes_{\dbZ} \dbZ_{(p)}}$ that fixes a family of endomorphisms of $L\otimes_{\dbZ} \dbZ_{(p)}$ (it seems to us that this first appears in [Zi]; see also [LR]). If $(G,X)$ is of $D$ type and $p$ is odd, then Kottwitz showed that $G_{\dbZ_{(p)}}$ is the connected component of the subgroup scheme of $\pmb{GL}_{L\otimes_{\dbZ} \dbZ_{(p)}}$ that fixes a family of endomorphisms of $L\otimes_{\dbZ} \dbZ_{(p)}$ (see [Ko]).

\bigskip
\noindent
{\bf 2.2. Hodge type.} We assume that there exists no reductive subgroup $G_1$ of $\pmb{GSp}(W,\psi)$ such that $G\subsetneq G_1\subseteq \pmb{GSp}(W,\psi)$ and $G_1^{\der}=G^{\der}$ (i.e., we assume that $(G,X)$ is saturated in $(\pmb{GSp}(W,\psi),S)$ in the sense of [Va1, Def. 4.3.1]). In [Va1, Prop 4.3.10 (a)] it is shown that there is $N(f)\in\dbN$ such  that if $p$ does not divide $N(f)$, then there exists a family of tensors $(v_{\alpha})_{\alpha\in\scrJ}$ of $\scrT(L\otimes_{\dbZ} \dbZ_{(p)})$ such that $G_{\dbZ_{(p)}}$ is the connected component of the subgroup scheme of $\pmb{GL}_{L\otimes_{\dbZ} \dbZ_{(p)}}$ that fixes $v_{\alpha}$ for all $\alpha\in\scrJ$. Each $v_{\alpha}$ is either an endomorphism of $L\otimes_{\dbZ} \dbZ_{(p)}$ or an element of $(L\otimes_{\dbZ} L\otimes_{\dbZ} L^*\otimes_{\dbZ} L^*)\otimes_{\dbZ} \dbZ_{(p)}$. The elements of $(L\otimes_{\dbZ} L\otimes_{\dbZ} L^*\otimes_{\dbZ} L^*)\otimes_{\dbZ} \dbZ_{(p)}$ are defined by perfect trace forms and involve Killing forms and Casimir elements. Therefore $N(f)$ is effectively computable.

If $(G,X)$ is not saturated in $(\pmb{GSp}(W,\psi),S)$, then the same holds provided we allow the tensors to be in $\scrT(L\otimes_{\dbZ} \dbZ_{(p)})$ (and not only in $(L\otimes_{\dbZ} L^*)\otimes_{\dbZ} \dbZ_{(p)}\cup (L\otimes_{\dbZ} L\otimes_{\dbZ} L^*\otimes_{\dbZ} L^*)\otimes_{\dbZ} \dbZ_{(p)}$). This is an easy exercise as $G_1/G$ is a torus and as tori over either $\dbQ$ or $\dbZ_{(p)}$ are linearly reductive (i.e., their representations are completely reducible).  

\bigskip\noindent
{\bf 2.3. General case.} 

\medskip\noindent
{\bf 2.3.1. Pseudo-Claim ([Ki3, p. 6]).} {\it There exists a family of tensors $(v_{\alpha})_{\alpha\in\scrJ}$ of $\scrT(L\otimes_{\dbZ} \dbZ_{(p)})$ such that $G_{\dbZ_{(p)}}$ is the subgroup of $\pmb{GL}_{L\otimes_{\dbZ} \dbZ_{(p)}}$ that fixes $v_{\alpha}$ for all $\alpha\in\scrJ$.}

\medskip\noindent
{\bf 2.3.2. Incomplete claim ([Ki4, Prop. 1.3.2]).} {\it Let $V$ be a discrete valuation ring of mixed characteristic $(0,p)$ and let $k$ be the residue field of $V$. Let $M$ be a free $V$-module of finite rank. Let $\scrG$ be an affine, flat, closed subgroup scheme of $\pmb{GL}_M$ such that the generic fibre of $\scrG$ is a reductive group scheme. Then $\scrG$ is the subgroup scheme of $\pmb{GL}_M$ that fixes a family of tensors of $\scrS(M)$.}

\medskip\noindent
{\bf 2.3.3. Remarks.} {\bf (a)}  A ``proof'' of ``Claim 2.3.2'' was made available in [Ki4, Prop. 1.3.2]. The proof of loc. cit. is incomplete as it relies on the the following false claim (see [Ki4, p. 10, l. 17]): ``each short exact sequence $0\to M_1\to M_2\to M_3\to 0$ of $\scrG$-modules with $\scrM_3$ as a trivial $\scrG$-module, splits as it splits after inverting $p$''. It is trivial to check that this ``claim'' is indeed false.

{\bf (b)} ``Claim 2.3.1'' was announced by Kisin around Jan. 2007 (see also [Ki3, p. 6]) but a proof of it never appeared. The approach to prove ``Claim 2.3.2'' used in [De3] and [Ki4] can not be used towards proving ``Claim 2.3.1'' .  

{\bf (c)} We are grateful to G. Pappas for pointing out to us that [Ki3] and [Ki4] contain different claims (see 2.3.1 versus 2.3.2).

{\bf (d)} One can correct ``Claim 2.3.2'' as follows (to be compared with [De3, Prop. 3.1 c)]). Suppose the homomorphism $\Hom(\pmb{GL}_M,\dbG_m)\to\Hom(\scrG,\dbG_m)$ has finite cokernel. Then as in [De3, Prop. 3.1 c)], one gets that in such a case the Claim 2.3.2 does hold. As the homomorphism $\Hom(\pmb{GL}_{L\otimes_{\dbZ} \dbZ_{(p)}},\dbG_m)\to\Hom(G_{\dbZ_{(p)}},\dbG_m)$ does have finite cokernel one gets:

\medskip\noindent
{\bf 2.3.4. Corollary.} {\it There exists a family of tensors $(v_{\alpha})_{\alpha\in\scrJ}$ of $\scrS(L\otimes_{\dbZ} \dbZ_{(p)})$ such that $G_{\dbZ_{(p)}}$ is the subgroup of $\pmb{GL}_{L\otimes_{\dbZ} \dbZ_{(p)}}$ that fixes $v_{\alpha}$ for all $\alpha\in\scrJ$.}

\bigskip\noindent
{\bf 2.4. Well positioned families of tensors.} Based on the limitations of Subsection 2.2 and of the false claim of Remark 2.3.3 (a), in [Va1] we developed a very general theory of well positioned families of tensors. Related to the context of Section 1, one has the following general and practical result (cf. [Va1, Prop. 4.3.10 (b)]).

\medskip\noindent
{\bf 2.4.1. Proposition.} {\it We assume that $p$ is odd and that the  Killing  form $\scrK$ on  $\Lie(G_{\dbZ_{(p)}}^{\der})$ and the form $\scrT$ on  $\Lie(G_{\dbZ_{(p)}}^{\der})$ induced (via restriction) by the trace form on $\text{End}(L\otimes_{\dbZ}\dbZ_{(p)})$ are both perfect. Let $(v_{\alpha})_{\alpha\in\scrJ_4}$ be the family of all elements of $\End(L\otimes_{\dbZ} \dbZ_{(p)})\cup (L\otimes_{\dbZ} L\otimes_{\dbZ} L^*\otimes_{\dbZ} L^*)\otimes_{\dbZ} \dbZ_{(p)}$ fixed by $G_{\dbZ_{(p)}}$. Let $\triangle$ be a faithfully flat $\dbZ_{(p)}$-algebra which is an integral domain. Let $M$ be a free $\triangle$-submodule of $L\otimes_{\dbZ} \triangle[{1\over p}]$ such that $M[{1\over p}]=L\otimes_{\dbZ} \triangle[{1\over p}]$ and we have $v_{\alpha}\in\scrT(M)$ for all $\alpha\in\scrJ_4$. Then the schematic  closure of $G_{\triangle[{1\over p}]}$ in $\pmb{GL}_M$ is a reductive group scheme over $\triangle$.}

\newpage
%\bigskip\smallskip
\noindent
{\boldsectionfont 3. Notations}

\medskip
Let $k$ be an algebraic closure of the field $\dbF_p$ with $p$ elements. Let $W(k)$ be the ring of Witt vectors with coefficients in $k$. Let $B(k):=W(k)[{1\over p}]$. Let $y:\Spec(k)\to\scrN$ be a point. Let $z:\Spec(V)\to\scrN$ be a lift of $y$, where $V$ is a discrete valuation ring that is a finite extension of $W(k)$. Let $e\in\dbN$ be the index of ramification of $V$. 

Let $\pi$ be a uniformizer of $V$. Let $R:=W(k)[[x]]$. Let $s:R\twoheadrightarrow V$ be a $W(k)$-epimorphism that maps $x$ to $\pi$. If $f_e$ is the Eisenstein polynomial of degree $e$ that has coefficients in $W(k)$ and that has $\pi$ as a root, then $\Ker(s)$ is generated by $f_e$. Let $S_e$ be the $R$-subalgebra of $B(k)[[x]]$ generated by all ${{x^{en}}\over {n!}}$ with $n\in\dbN\cup\{0\}$; it is the divided power hull of any one of the ideals $(x^e)$, $(f_e)$, or $(p,x^e)=(p,f_e)$ of $R$. Let $J_e$ (resp. $K_e$) be the ideal of $S_e$ generated by all ${{f_e^{en}}\over {n!}}$ (resp. by all ${{x^{en}}\over {n!}}$) with $n\in\dbN$. Let $R_e$ be the $p$-adic completion of $S_e$. Let $\Rtil_e$ be the completion of $S_e$ with respect to the decreasing filtration given by its ideals $K_e^{[n]}$, $n\in\dbN\cup\{0\}$. Thus $\Rtil_e=\text{proj}.\text{lim}._{n\in\dbN} S_e/K_e^{[n]}$. We recall that $K_e^{[0]}:=S_e$ and that for $n\Ge 1$ the ideal $K_e^{[n]}$ of $S_e$ is generated by all products $\frac{\delta_1^{a_1}}{a_1!}\cdots\frac{\delta_m^{a_m}}{a_m!}$ with $\delta_1,\ldots,\delta_m\in K_e$ and with $m$, $a_1,\ldots,a_m\in\dbN\cup\{0\}$ such that we have $a_1+\cdots+a_m\ge n$. If $p>2$, then $\Rtil_e$ is also the completion of $S_e$ with respect to its decreasing filtration $(J_e^{[n]})_{n\in\dbN\cup\{0\}}$. Thus for $p\ge 3$ we have as well a $W(k)$-epimorphism $\tilde e_V:\Rtil_e\twoheadrightarrow V$ that takes $x$ to $\pi$ and for $p=2$ we have a $W(k)$-epimorphism $\tilde e_V:\Rtil_e\twoheadrightarrow V/pV$ that takes $x$ to $\pi$ modulo $p$.

The $W(k)$-algebra $\Rtil_e$ (resp. $R_e$) consists
of formal power series $\Sigma_{n\Ge 0} a_nx^n$ such that the sequence $([{n\over e}]!a_n)_{n\in\dbN\cap\{0\}}$ is formed by elements of $W(k)$ (resp. is formed by elements of $W(k)$ and converges to $0$). Let $\Phi_k$ be the Frobenius lift of $R$, $S_e$, $R_e$, or $\Rtil_e$ that is compatible with $\sigma$ and such that
$\Phi_k(x)=x^p$. One has $R_k\subseteq \Rtil_e$ and $\Phi_k(\Rtil_e)\subseteq R_e$.

Let $O_R$ be the unique local ring of $R$ that is a discrete valuation ring of mixed characteristic $(0,p)$ and residue field $k((x))$. Let $\scrO$ be the completion of $O_R$. Let $k_1$ be an algebraic closure of $k((x))$. Let
$$\Spec(W(k_1))\to \Spec(R)$$ 
be the Teichm\"uller lift with respect to $\Phi_k$; under it $W(k_1)$ gets naturally the structure of a $*$-algebra, where $*\in\{R,O_R,\scrO\}$. 

To $z:\Spec(V)\to\scrN$ corresponds a principally polarized abelian scheme $(A,\lambda_A)$ over $V$ whose generic fibre is endowed with a family $(w_{\alpha})_{\alpha\in\scrJ}$ of Hodge cycles. Let $(N,\phi_N,\psi_N,\nabla_N)$ be the evaluation of $\dbD(((A,\lambda_A)[p^{\infty}])_{V/pV})$ at the thickening associated naturally to the closed embedding $\Spec(V/pV)\hookrightarrow \Spec(R_e)$. Thus $N$ is a free $R_e$-module of rank $2d$, $\phi_N:N\to N$ is a $\Phi_k$-linear endomorphism, $\psi_N$ is  a perfect alternating form on $N$ which is a principal quasi-polarization of $(N,\phi_N)$, and $\nabla_N:N\to N\otimes_{R_e} R_edX$ is a connection on $N$ with respect to which $\phi_N$ is horizontal i.e., we have $\nabla_N\circ\phi_N=(\phi_N\otimes d\Phi_k)\circ\nabla_N$. The connection $\nabla_N$ is integrable and nilpotent modulo $p$.

We have a functorial identification $H^1_{dR}(A_V/V)=N\otimes_{R_e} V$. Let $F^1_V$ be the Hodge filtration of $H^1_{dR}(A_V/V)$ defined by $A_V$. 

Let $(M_0,\phi_0,(t_{0,\alpha})_{\alpha\in\scrJ},\psi_0)$ be the tensorization of $(N,\phi_N,(t_{z,\alpha})_{\alpha\in\scrJ},\psi_N)$ via the epimorphism $\xi:R_e\twoheadrightarrow W(k)$ that maps $x$ to $0$. It is known that there exists isomorphisms $(M_0,(t_{0,\alpha})_{\alpha\in\scrJ})\arrowsim (L^*\otimes_{\dbZ} W(k),(v_{\alpha})_{\alpha\in\scrJ})$ and that there exists a unique isomorphism $(N_0[{1\over p}],(t_{z,\alpha})_{\alpha\in\scrJ})\arrowsim (M_0\otimes_{W(K)} R_e[{1\over p}],(t_{0,\alpha})_{\alpha\in\scrJ})$ which modulo $\Ker(\xi[{1\over p}])$ is the identity automorphism of $M_0[{1\over p}]$.

For $\alpha\in\scrJ$, the crystalline realization of $w_{\alpha}$ is a tensor $t_{z,\alpha}$ of $\scrT(N)[{1\over p}]$ which is fixed by $\phi_N$, which is annihilated by $\nabla_N$, and whose tensorization with $V[{1\over p}]$ belongs to the $F^0$-filtration of $\scrT(H^1_{dR}(A_V/V))[{1\over p}]$ defined by $F^1_V[{1\over p}]$. Let $\tilde G_{R_e}$ be the schematic closure in $\pmb{GL}_N$ of the reductive subgroup scheme of $\pmb{GL}_{N[{1\over p}]}$ that fixes $t_{z,\alpha}$ for all $\alpha\in\scrJ$ (this makes sense as there exists isomorphisms $(N[{1\over p}],(t_{z,\alpha})_{\alpha\in\scrJ})\arrowsim (L^*\otimes_{\dbZ} W(k),(v_{\alpha})_{\alpha\in\scrJ})$).

Let $\eta$ be the field of fractions of $R$. Let $U:=\Spec(R)\setminus\Spec(k)$ be the open subscheme of $\Spec(R)$ which is the complement in $\Spec(R)$ of its closed point $\Spec(k)$. 

\bigskip\smallskip
\noindent
{\boldsectionfont 4. Method I}

\medskip
Method I is due to Faltings and Vasiu, works for $p\ge 5$, and uses the following tool.

\bigskip\noindent
{\bf 4.1. Faltings theorem (see [Fa, Cor. 9]).} {\it Let $\alpha\in\scrJ$ be such that we have $v_{\alpha}\in (\cup_{i=1}^{p-2} L^{\otimes i}\otimes_{\dbZ} L^{* \otimes i}\otimes_{\dbZ} \dbZ_{(p)})$. Then we have also $t_{z,\alpha}\in (\cup_{i=1}^{p-2} N^{\otimes i}\otimes_{R_e} N^{* \otimes i})$.}

\bigskip\noindent
{\bf 4.2. Step 1 (the reductiveness part).}  The first  step is to show that, under some conditions on the closed embedding homomorphism $G_{\dbZ_{(p)}}\hookrightarrow \pmb{GSp}(L\otimes_{\dbZ}\dbZ_{(p)},\psi)$ and under the assumption that $p\ge 5$, $\tilde G_{R_e}$ is a reductive subgroup scheme of $\pmb{GL}_N$. 
See [Va1, Subsect. 5.2] for more details and see [Va1, (5.2.12)] for the fact that the reductive group scheme $\tilde G_{R_e}$ is isomorphic to $G_{\dbZ_{(p)}}\times_{\dbZ_{(p)}} R_e$.

Here is a basic example of the type of conditions one requires.
By combining Theorem 4.1 with Proposition 2.4 one gets. 

\smallskip\noindent
{\bf 4.2.1. Corollary.} {\it We assume that $p\ge 5$ and that the  Killing  form $\scrK$ on  $\Lie(G_{\dbZ_{(p)}}^{\der})$ and the form $\scrT$ on  $\Lie(G_{\dbZ_{(p)}}^{\der})$ induced (via restriction) by the trace form on $\text{End}(L\otimes_{\dbZ}\dbZ_{(p)})$ are both perfect. Then $\tilde G_{R_e}$ is a reductive group scheme over $R_e$.}

\bigskip\noindent
{\bf 4.3. Step 2 (lift of the filtration).} The second  step shows that we can lift $F^1_V$ to a direct summand $F^1_{R_e}$ of $N$ in such a way that $\psi_N(F^1_{R_e},F^1_{R_e})=0$ and that for each element $\alpha\in\scrJ$ the tensor $t_{z,\alpha}$ belongs to the $F^0$-filtration of $\scrT(N)[{1\over p}]$ defined by  $F^1_{R_e}[{1\over p}]$. The essence of this second  step is the classical theory of infinitesimal liftings of cocharacters of smooth group schemes (see [DG, Exp. IX]). In other words, $F^1_V$ is defined by a suitable cocharacter of $\tilde G_{R_e}\times_{R_e} V$ which lifts to a cocharacter of $\tilde G_{R_e}$ that defines $F^1_{R_e}$. See [Va1, Subsect. 5.3] for more details. Strictly speaking, loc. cit. is worked out over $\Rtil_e$ instead of over $R_e$; but as for each $n\in\dbN$ the $W(k)$-algebra $R_e/p^nR_e$ is the inductive limit of its local artinian $W(k)$-subalgebras, the arguments of loc. cit. work over $R_e$ itself instead of over $\Rtil_e$. Due to the existence of $F^1_{R_e}$, for $p\ge 3$ the composite morphism $z:\Spec(V)\to \scrN\to\scrM$ lifts to a morphism $\tilde z:\Spec(R_e)\to \scrM$ (cf. classical deformation theories of Grothendieck--Messing and Serre--Tate). The first versions of [Fa] and [Va1] used $R_e$ and $\Rtil_e$ (respectively). The final version of [Fa] was worked out over $\Rtil_e$.  

\bigskip\noindent
{\bf 4.4. Step 3 (parallel transport of Hodge cycles).} To $\tilde z$ corresponds an abelian scheme over $R_e$ whose pull back to $R_e[{1\over p}]$ is endowed with a family of Hodge cycles whose crystalline realizations are $(t_{z,\alpha})_{\alpha\in\scrJ}$. This is a result of Faltings whose essence is presented in [Va1, Rm. 4.1.5] and whose detailed proof is presented in [Va8, Prop. 3.4.1] (one can embed $R_e$ into $\dbC[[x]]$ and one can check the statement over rings of formal power series with coefficients in $\dbC$; this case is the same as [Va8, Prop. 3.4.1]). This implies that the morphism $\tilde z:\Spec(R_e)\to \scrM$ factors naturally through $\scrN^{\text{cl}}$. If $R_e^{\text{n}}$ is the normalization of $R_e$, then one gets a morphism $\tilde z^{\text{n}}:\Spec(R_e^{\text{n}})\to \scrN$ that lifts $z$. Using this and the fact that we have a $W(k)$-epimorphism $R_e^{\text{n}}\twoheadrightarrow W(k)$ that maps $x$ to $0$, we conclude that there exists a lift $z_0:\Spec(W(k))\to\scrN$ of $y$.

\bigskip\noindent
{\bf 4.5. Step 4 (deformation theory).} The fourth step uses the lift $z_0:\Spec(W(k)))\to\scrN$ of $y$ and Faltings deformation theory (see [Fa, Sect. 7]) to show that $\scrN$ is formally smooth over $\dbZ_{(p)}$ at its $\dbF$-valued point defined by $y$. See [Va1, Subsect. 5.4] for more details (which require the form of Step 3 over rings of formal power series over $W(k)$).

\bigskip\noindent
{\bf 4.6. Step 5 (well positioned families of tensors).} The fifth  step shows that for $p>3$ the mentioned conditions on the closed embedding homomorphism $G_{\dbZ_{(p)}}\hookrightarrow \pmb{GSp}(L\otimes_{\dbZ}\dbZ_{(p)},\psi)$ always hold, {\it provided} we replace $f:(G,X)\hookrightarrow (\pmb{GSp}(W,\psi),S)$ by a suitable other injective map $f_1:(G_1,X_1)\hookrightarrow (\pmb{GSp}(W_1,\psi_1),S_1)$ with the property that $(G^{\ad},X^{\ad})=(G^{\ad}_1,X^{\ad}_1)$ and that we have an isogeny $G_1^{\der}\to G^{\der}$. See [Va1, Subsects. 6.5 and 6.6] for more details. Using this, from the fact that Theorem 1.1 holds for the analogue $\scrN_1$ of $\scrN$ (obtained working with $f_1$ and a prime $v_1$ of $E(G_1,X_1)$ that divides the same prime of $E(G_1^{\ad},X_1^{\ad})=E(G^{\ad},X^{\ad})$ as $v$), one gets directly that Theorem 1.1 holds for $\scrN$ itself (see [Va1, Subsubsects. 3.2.12, 6.1, and 6.2]. 

\medskip\noindent
{\bf 4.6.1. Example.} We assume that $p\ge 5$, that $G$ is a $\pmb{GSpin}(2,2n-1)$ group, and that $d$ is a power of $2$. We have $E(G,X)=\dbQ$ and $O=\dbZ_{(p)}$. 

If $p$ does not divide $2n-1$, then the Killing  form $\scrK$ on  $\Lie(G_{\dbZ_{(p)}}^{\der})$ and the form $\scrT$ on  $\Lie(G_{\dbZ_{(p)}}^{\der})$ induced (via restriction) by the trace form on $\text{End}(L\otimes_{\dbZ}\dbZ_{(p)})$ are both perfect. Therefore Corollary 4.2.1 applies and thus Theorem 1.1 holds. 

If $p$ divides $2n-1$, then it does not divide $2n+1$ and we can assume that $f$ is a composite injective map $(G,X)\hookrightarrow (G_1,X_1)\hookrightarrow (\pmb{GSp}(W,\psi),S)$ such that $G_1$ is a $\pmb{GSpin}(2,2n+1)$ group and the schematic closure $G_{1,\dbZ_{(p)}}$ of $G_1$ in $\pmb{GL}_{L\otimes_{\dbZ} \dbZ_{(p)}}$ is a reductive group scheme. Let $H_{1,p}:=G_1(\dbQ_p)\cap K_p$; it is a hyperspecial subgroup of $G_1(\dbQ_p)$. As $p$ does not divide $2n+1$, the normalization $\scrN_1$ of the schematic closure $\scrN_1^{\text{cl}}$ of $\Sh(G_1,X_1)/H_{1,p}$ in $\scrM$ is a regular scheme which is formally smooth over $O$ (cf. previous paragraph); moreover the analogue $\tilde G_{1,R_e}$ of $\tilde G_{R_e}$ is a reductive group scheme over $R_e$ which is a closed subgroup scheme of $\pmb{GL}_N$. There exists a rank $1$ torus $T_{\dbZ_{(p)}}$ of $G_{1,\dbZ_{(p)}}$ such that the centralizer of $T_{\dbZ_{(p)}}$ in $G_{1,\dbZ_{(p)}}$ is a reductive group scheme $G_{1/2,\dbZ_{(p)}}$ that contains $G_{\dbZ_{(p)}}$ in such a way that $G^{\der}_{1/2,\dbZ_{(p)}}=G^{\der}_{\dbZ_{(p)}}$. Let $\scrN_{1/2}$ be defined similarly to $\scrN_1$ but working with $H_{1/2,p}:=G_{1/2}(\dbQ_p)\cap K_p$ and $\Sh(G_{1/2},X_{1/2})/H_{1/2,p}$. The injective map  $(G_{1/2},X_{1/2})\hookrightarrow (G_1,X_1)$ is a {\it relative PEL embedding} and therefore (as in [Zi] and [Ko]) one argues that $\scrN_{1/2}$ is a regular scheme which is formally smooth over $O=\dbZ_{(p)}$ and which is a closed subscheme of $\scrN_1$; moreover the analogue $\tilde G_{1/2,R_e}$ of $\tilde G_{R_e}$ is a reductive group scheme over $R_e$ which is a closed subgroup scheme of $\tilde G_{1,R_e}$ and thus of $\pmb{GL}_N$. As $\scrN$ is an open subscheme of $\scrN_{1/2}$ we conclude that $\scrN$ is a regular scheme which is formally smooth over $O=\dbZ_{(p)}$ and that $\tilde G_{R_e}$ is a reductive group scheme over $R_e$.

\bigskip\noindent
{\bf 4.7. Remarks.} {\bf (a)} From the smoothness of $\scrN$ and Faltings deformation theory one gets as a corollary that $t_{0,\alpha}\in\scrT(M_0)[{1\over p}]$ does not depend on the lift $z:\Spec(V)\to\scrN$ of $y:\Spec(\dbF)\to\scrN$.

{\bf (b)} If $f:(G,X)\to (\pmb{GSp}(W,\psi),S)$ is a PEL type embedding (with the $D$ case excluded for $p=2$), then the formal smoothness of $\scrN^{\text{cl}}$ over $O$ at the point $y$ follows from a direct computation of the tangent space of $\scrN$ at $y$. The general principle that {\it reductiveness of $\tilde G_{R_e}$} implies the {\it formal smoothness of $\scrN$ at $y$} appears first in [Va0] and thus in [Va1] and all its preliminary versions. 

\bigskip\noindent
{\bf 4.8. Comments.} The main merits of Method I are that:

\smallskip
{\bf  (a)} it provides a general theory of well positioned family of tensors (which gives a lot more information than just the reductiveness of $\tilde G_{R_e}$);

\smallskip
{\bf (b)} it introduces the fundamental concept of a {\it relative PEL embedding} (see [Va1, Subsubsects. 4.3.16, 4.3.13, 4.3.11]);

\smallskip
{\bf (c)} it works uniformly for all primes $p\ge 5$ and it is able to avoid all the complications of ``Claims 2.3.1 and 2.3.2'' (despite the relatively weak tool provided by Theorem 4.1).

\medskip
Among the weak points of Method I, we mention that: {\bf (i)}  it involves cases (on Shimura types and on the prime $p$, like Example 4.6.1 pertain to the $B_n$ type and one had to consider two cases), {\bf (ii)} (due to the complicated nature of $R_e$) it requires some more expertise on reductive group schemes, and {\bf (iii)} it does not work for $p\in\{2,3\}$. Therefore we felt the need to simplify Method I and to find a new method that works as well for $p\in\{2,3\}$. A key idea of Method II is to replace the role of $R_e$ but that one of the regular, local ring $R=W(k)[[x]]$ (see Section 5 below).

\bigskip\noindent
{\bf 4.9. Moonen's report.} There has been an erroneous report on [Va1], cf. [Mo]. On one side [Mo] contains many gaps, errors, and mistakes and on the other side [Mo] does not provide any concrete or serious issue with [Va1]. See [Va3] for errata to [Va1]. Except a gap in Faltings work that was corrected in [Va2] and a minor confusing for the $A_n$ type case with $p$ dividing $n+1$ that was corrected in [Va3]${}^1$ $\vfootnote{1}{The confusion pertained only to the passage from Hodge type to the abelian type and not to cases related to Theorem 1.1.}$, there has never been any other trouble with the Method I presented in [Va1]. Below we present only those gaps, errors, and mistakes in [Mo] which had been pointed out by us to Moonen${}^2$ $\vfootnote{2}{They were pointed out before [Mo] was accepted for publication. Moonen informed us that he prefers to let (quoting)  ``the specialists to decide''.}$  and those who had been of concern to us since many years.${}^3$ $\vfootnote{3}{Moonen informed us that the only thing he cares about is the mathematical truth.}$ 

\medskip
{\bf (a)} Proof of [Mo, Prop. 3.22] is wrong. It ``proves'' that any
cyclic Galois cover $K_1$ of the field of fractions $K$ of a discrete valuation ring of mixed
characteristic, is automatically unramified. The error
consists in the fact that if an element of the Galois group acts trivially
on the residue field of the ring of integers $O_1$ of $K_1$, it does not
necessarily act trivially on $O_1/pO_1$. Any attempt of handling [Va1, Thm.
6.1.2*] only at the level of discrete valuation rings (i.e., as in [Mo]) is pointless (see [Va1, Subsect. 6.2.7]).

{\bf (b)} The descent arguments of [Mo, Prop. 3.10] are not OK.
The counterexample of [BLR, Ch. 6, Sect. 6.7] makes perfect sense in mixed
characteristic.

{\bf (c)} The proof of Faltings Lemma in [Mo, Lemma 3.6] is wrong. It ends up
that (quoting) ``At some points one furthermore needs arguments similar to the above, i.e., taking sections over an extension ... We leave it to the reader to check the details" [in the
book of Faltings and Chai]. But as pointed out in [Va2], there has been a gap in Faltings work (including his book
with Chai) and that consisted exactly in the part of taking sections (Faltings erroneous argument was reproduced in the third paragraph of [Va1, Step B, Subsubsect. 3.2.17]). The gap was corrected only in [Va2] using two proofs: one based on smooth toroidal compactifications and one based on de Jong extension theorem of [dJ], cf. [Va2, Prop. 4.1 and Rm. 4.2].

{\bf (d)}  Obviously, (a) to (c) invalidate many other parts of [Mo], like
[Mo, Prop. 3.18 and Cor. 3.23]. In fact they invalidate
essentially all of [Mo, Section 3] that pertains to integral models. The
very few parts of [Mo, Section 3] that pertain to integral models and that
are correct, are the only ones that followed entirely [Va1] and its
preliminary versions. 

{\bf (e)} Based on (a) to (d), the claim of [Mo, Rm. 3.24] of a ``very
different" presentation from the earlier versions of [Va1], is
meaningless.

\bigskip\smallskip
\noindent
{\boldsectionfont 5. Method II}

\medskip
Let $\scrN^{\text{s}}$ be the formally smooth locus of $\scrN$ over $O$. Method II is due to Vasiu, works for all primes $p$, and uses the following tool.

\bigskip\noindent
{\bf 5.1. de Jong's Theorem (see [dJ, Thm. 1.1])}. {\it The natural functor from the category of $F$-crystals over $\Spec(k[[x]])$ to the category of $F$-crystals over $\Spec(k((x)))$ is fully faithful.}

\medskip
The main idea of Method II is to work out the reductiveness part not over the complicated rings $R_e$ but over the regular local ring $R$ of dimension $2$. Theorem 5.1 is exactly what one needs in order to accomplish the replacement of $R_e$ by $R$.  

\bigskip\noindent
{\bf 5.2. Theorem (Step 1, a motivic conjecture of Milne proved in [Va5, Thm. 1.3]).} {\it We assume that either $p>2$ or $G$ is a torus. We also assume that $V=W(k)$. Then there exists an isomorphism $(M_0,(t_{0,\alpha})_{\alpha\in\scrJ})\arrowsim (L^*\otimes_{\dbZ} W(k),(v_{\alpha})_{\alpha\in\scrJ})$.}

\bigskip\noindent
{\bf 5.3. Step 2 (deformation theory).} Faltings deformation theory shows that $\scrN$ is smooth at $y$ provided we can choose $V=W(k)$ (this holds independently of the first step and without assuming that $G_{\dbZ_{(p)}}$ is a reductive subgroup scheme of $\pmb{GL}_{L\otimes_{\dbZ} \dbZ_{(p)}}$, cf. [Va8, Thm. 1.5 (a)]). Using maximal tori of $G_{\dbZ_{(p)}}$ whose pull backs to $\dbR$ have $\dbR$-rank $1$, one can easily get that $\scrN$ has $W(k)$-valued points and thus that $\scrN^{\text{s}}$ has a non-empty intersection with $\scrN_{k(v)}$. 

\bigskip\noindent
{\bf 5.4. Step 3 (the reductiveness and the lift of filtration part).} The second step shows that $\scrN^{\text{s}}$ intersects $\scrN_{k(v)}$ in  an open closed subscheme of $\scrN_{k(v)}$. To check this one only has to show that for each commutative diagram of the following type  
$$
\CD
\Spec(k) @>{}>> \Spec(k[[x]]) @<{}<<
\Spec(k((x))) \\
@VV{y}V @VV{q}V @VV{q_{k((x))}}V \\
\scrN @<{}<<  \scrN @<{}<< \scrN^{\text{s}},
\endCD
$$ 
the morphism $y:\Spec(k)\to\scrN$ (or $q$) factors through the open subscheme $\scrN^{\text{s}}$ of $\scrN$.

We consider the principally quasi-polarized $F$-crystal
$$(N_0,\phi_{N_0},\nabla_0,\psi_{N_0})$$ 
over $k[[x]]$ of the principally polarized abelian scheme over $k[[x]]$ associated to $q$. Thus $N_0$ is a free $R$-module of rank $2d$, $\phi_{N_0}$ is a $\Phi_k$-linear endomorphism of $N_0$ such that $p$ annihilates $N_0/R\text{Im}(\phi_{N_0})$, $\psi_{N_0}$ is  a perfect alternating form on $N_0$ which is a principal quasi-polarization of $(N_0,\phi_{N_0})$, and $\nabla_0$ is an integrable and nilpotent modulo $p$ connection on $N_0$ such that we have $\nabla_0\circ\phi_{N_0}=(\phi_{N_0}\otimes d\Phi_k)\circ\nabla_0$.

As the $O$-scheme $\scrN^{\text{s}}$ is formally smooth, there exists a lift $z_1:\Spec(\scrO)\to\scrN^{\text{s}}$ of the morphism $q_{k((x))}:\Spec(k((x)))\to\scrN^{\text{s}}$ defined naturally by $q_{k((x))}$ and denoted in the same way. Let 
$(A_1,\lambda_{A_1},(w_{1,\alpha})_{\alpha\in\scrJ})$
be the principally polarized abelian scheme over $\scrO$ endowed with a family of Hodge cycles that is associated naturally to $z_1$. Let $t_{1,\alpha}$ be the de Rham realization of $w_{1,\alpha}$. We identify canonically $N_0\otimes_R \scrO=H^1_{\text{dR}}(A_1/\scrO)$ and thus we can view each $t_{1,\alpha}$ as a tensor of $\scrT(N_0\otimes_R \scrO)[{1\over p}]$. Based on Theorem 5.1 one easily gets that (see [Va8, Part I, Prop. 5.1.1]):

\medskip
{\bf (*)} in fact we have $t_{1,\alpha}\in\scrT(N_0)[{1\over p}]$ for all $\alpha\in\scrJ$.

\medskip
Based on (*), one can speak about the reductive subgroup $\tilde G_{\eta}$ of $\pmb{GL}_{N_0\otimes_R \eta}$ that fixes each  $t_{1,\alpha}$ with $\alpha\in\scrJ$. Let $\tilde G_R$ be the schematic closure of $\tilde G_{\eta}$ in $\pmb{GL}_{N_0}$. From Theorem 5.2 we get that (see [Va8, Part I, Thm. 5.2]):

\medskip
{\bf (**)} there exists an isomorphism $(N_0\otimes_R W(k_1),(t_{1,\alpha})_{\alpha\in\scrJ})\arrowsim (L^*\otimes_{\dbZ} W(k_1),(v_{\alpha})_{\alpha\in\scrJ})$ and thus both $\tilde G_R\times_R W(k_1)$ and $\tilde G_R\times_R \scrO$ are reductive group schemes. 

\medskip
From this and the crystalline reductiveness principle of [Va6, Thm. 6.3] we get that:

\medskip
{\bf (***)} $\tilde G_R$ is a reductive group scheme over $R$ (we recall that $R=W(k)[[x]]$).

\medskip
We recall that the crystalline reductiveness principle relies on the fact that each reductive group scheme over $U$ extends uniquely to a reductive group scheme over $\Spec(R)$ (see [Va6, Thm. 1.1 (c)] whose proof relies on [CTS]). Based on (***), as in Subsection 4.3 one argues that one can lift the kernel of $\phi_{N_0}$ modulo $p$ to a direct summand $F^1_R$ of $N_0$ in such way that each tensor $t_{1,\alpha}$ belongs to the $F^0$-filtration of $\scrT(N_0)[{1\over p}]$ defined by $F^1_R$. Using this and Faltings deformation theory one gets that indeed $q$ factors through $\scrN^{\text{s}}$ (see [Va8, Part I] for $p\ge 3$ and for some exceptional cases with $p=2$ and see [Va8, Part II] for the general case with $p=2$).

\bigskip\noindent
{\bf 5.6. Step 4 (the three cases).} Each simple, adjoint Shimura pair $(G_0,X_0)$ of abelian type belongs to precisely one of the following three classes:

\medskip
{\bf Type (I)} $(G_0,X_0)$ is of $A_n$, $C_n$, or $D_n^{\dbH}$ type and $G_{0,\dbR}$ has no simple, compact factor.

\smallskip
{\bf Type (II)} The group $G_{0,\dbR}$ has a simple, compact factor.

\smallskip
{\bf Type (III)} $(G_0,X_0)$ is of $B_n$ (with $n\ge 3$) or $D_n^{\dbR}$ (with $n\ge 4$) type and $G_{0,\dbR}$ has no simple, compact factor.

\medskip
If $(G^{\ad},X^{\ad})=(G_0,X_0)$, then for each one of the mentioned three classes one can choose $f:(G,X)\hookrightarrow (\pmb{GSp}(W,\psi),S)$ such that for all primes $v$ dividing $p$ the scheme $\scrN$ is regular and formally smooth over $O$ (i.e., the Theorem 1.1 holds). Concretely, we have the following three cases.

\medskip
{\bf (i)} If $(G_0,X_0)$ is of Type I, then we can assume that $f:(G,X)\hookrightarrow (\pmb{GSp}(W,\psi),S)$ is a PEL type embedding and that the centralizer of $G_{\dbZ_{(p)}}$ in $\pmb{GL}_{L\otimes_{\dbZ} \dbZ_{(p)}}$ is the reductive group scheme of invertible elements of a semisimple $\dbZ_{(p)}$-subalgebra of $\End(L\otimes_{\dbZ} \dbZ_{(p)})$. See [Va3, Prop. 3.2] for the $A_n$ type case and see [Va8, Part II] for the $C_n$ and $D_n^{\dbH}$ types (the arguments are the same as in [Va3, Prop. 3.2]). Except when $p=2$ and $(G_0,X_0)$ is of $D_n^{\dbH}$ type, it is well known that Theorem 1.1 holds (cf. [Zi] and [Ko]). The case when $p=2$ and $(G_0,X_0)$ is of $D_n^{\dbH}$ type is proved in [Va7].

\medskip
{\bf (ii)} If $(G_0,X_0)$ is of Type II, then for any choice $f:(G,X)\hookrightarrow (\pmb{GSp}(W,\psi),S)$ one knows that $\scrN$ is a pro-\'etale cover of a projective $O$-scheme $\scrN/H_0$ (cf. [Va4, Cor. 4.3]) and that the connected components of $\scrN$ are permuted transitively by $G(\dbA_f^{(p)})$ (cf. [Va1, Lemma 3.3.2]).  From this and the fact that $\scrN$ has $W(k)$-valued points (see Subsection 5.3) we get that each special fibre of $\scrN/H_0\times_O W(k)$ is connected and that the connected components of $\scrN_{k(v)}$ are permuted transitively by $G(\dbA_f^{(p)})$. Based on this and the fact that $\scrN^{\text{s}}/H_0\cap\scrN_{k(v)}/H_0$ is a $G(\dbA_f^{(p)})$-invariant non-empty open closed subscheme of $\scrN_{k(v)}/H_0$, we  easily conclude that $\scrN^{\text{s}}/H_0=\scrN/H_0$ and that $\scrN^{\text{cl}}=\scrN$. Thus Theorem 1.1 holds.

\medskip
{\bf (iii)} If $(G_0,X_0)$ is of Type III, then we can assume that $f:(G,X)\hookrightarrow (\pmb{GSp}(W,\psi),S)$ is such that $E(G,X)=\dbQ$ and that the ordinary locus of $\scrN_{k(v)}$ is Zariski dense in $\scrN_{k(v)}$.${}^4$ $\vfootnote{4}{This and the $p=2$ case of (i) with $(G_0,X_0)$ of $D_n^{\dbH}$ type and with $G_{\dbQ_2}$ split (see [Va7, Part I]), are the only places in Method II where one needs to use ramification rings like $V$ and $R_e$. If one knows that the connected components of $\scrN_{k(v)}$ are permuted transitively by $G(\dbA_f^{(p)})$ or if one is eager to use toroidal induction in order to prove this transitivity property, then one does not even have to mention $V$ and $R_e$ at all.}$ As each ordinary point of $\scrN$ belongs to $\scrN^{\text{s}}$ (cf. [No, Cor. 3.8]), one concludes (cf. Subsection 5.4) that $\scrN=\scrN^{\text{s}}$. Thus Theorem 1.1 holds. See [Va8, Part II] for details.

\medskip
The general case of Theorem 1.1 follows from the above three cases via the standard isogeny properties between Shimura varieties of Hodge type (see [Va1], [Va3], and [Va8] for details).

\bigskip\noindent
{\bf 5.7. Comments.} The main merits of Method II are that:

\smallskip
{\bf  (a)} it works for all primes $p$;

\smallskip
{\bf (b)} it works over $R$ and not over $R_e$ (and this replacement of $(V,R_e)$ by $(W(k),R)$ is almost complete, cf. footnote 4);

\smallskip
{\bf (c)} it can be easily applied to other classes of polarized varieties and (the greatest part of) it can be adapted to many situations in which $G_{\dbZ_{(p)}}$ is not a reductive group scheme over $\dbZ_{(p)}$.

\medskip
We would be inclined to think that it has no weak point. The only arguable think would be that, in some sense, it still involves cases (see Types I to III above). But if one really wants to handle the case $p=2$ as well, then (due to so many pathologies in characteristic $2$ and the desire to keep the things relatively simple) it seems unavoidable to have some cases to consider. 

\bigskip\smallskip
\noindent
{\boldsectionfont 6. Method III}

\medskip
Method III is in essence only a variation of the Methods I and II. Method III is due to Kisin, works for $p>2$ and for some exceptional cases with $p=2$, and uses the following tool.

\bigskip\noindent
{\bf 6.1. Kisin's Theorem (see [Ki1] and [Ki2]).} {\it If $p=2$ we assume that the $2$-rank of $A_k$ is $0$. Then the quadruple $(N,\phi_N,(t_{z,\alpha})_{\alpha\in\scrJ},\psi_N)$ is the extension via the $\sigma$-linear homomorphism $R\to R_e$ that maps $x$ to $x^p$ of a quadruple $(N_0,\phi_{N_0},(t_{z,\alpha})_{\alpha\in\scrJ},\psi_{N_0})$ over $R$, where $\phi_{N_0}:N_0\to N_0$ is a $\sigma$-linear endomorphism such that $N_0/R\text{Im}(\phi_{N_0})$ is annihilated by $f_e(x)$. Moreover, for $\alpha\in\scrJ$ we have $v_{\alpha}\in \scrT(L\otimes_{\dbZ} \dbZ_{(p)})$ if and only if $t_{z,\alpha}\in\scrT(N_0)$.} 

\medskip
Obviously Theorem 6.1 implies Theorem 4.1 and therefore it is a significantly better tool. We have the following analogue of the property 5.4 (*):

\medskip
{\bf (*)} we have $t_{z,\alpha}\in\scrT(N_0)[{1\over p}]$ for all $\alpha\in\scrJ$.

\medskip
Based on (*), one can speak about the reductive subgroup $\tilde G_{\eta}$ of $\pmb{GL}_{N_0\otimes_R \eta}$ that fixes each  $t_{z,\alpha}$ with $\alpha\in\scrJ$. Let $\tilde G_R$ be the schematic closure of $\tilde G_{\eta}$ in $\pmb{GL}_{N_0}$.

\bigskip\noindent
{\bf 6.2. Step 1 (Milne conjecture).} In [Fo] it is proved that the following analogue of the property 5.4 (**) holds:

\medskip
{\bf (**)} we have an isomorphism $(N_0\otimes_R W(k_1),(t_{1,\alpha})_{\alpha\in\scrJ})\arrowsim (L^*\otimes_{\dbZ} W(k_1),(v_{\alpha})_{\alpha\in\scrJ})$ and thus $\tilde G_R\times_R \scrO$ is a reductive group scheme. 

\bigskip\noindent
{\bf 6.3. Step 2 (the reductiveness part).} One has the following analogue of the property 5.4 (***). 

\medskip
{\bf (***)}  $\tilde G_R$ is a reductive group scheme over $R$ (we recall that $R=W(k)[[x]]$).

\medskip
Due to ``Claim 2.3.2'', [Ki4, Prop. 1.3.4] is of little use as it stands. Thus, in order to eliminate the (strong) hypothesis of [Ki4, Prop. 1.3.4] that pertains to a reductive groups scheme over $\dbZ_p$ defined by a family of tensors, the proof of [Ki4, Prop. 1.3.4] needs modifications. We present two proofs of (***) (and implicitly of [Ki4, Prop. 1.3.4]) that do not appeal to the mentioned (strong) hypothesis and which give the same result. 

{\bf Proof 1.} The proof of the crystalline reductiveness principle of [Va6, Thm. 6.3] applies entirely to give us that (***) holds (once 6.1 (*) holds, it is irrelevant in the proof of loc. cit. if $p$ or if $f_e(x)$ annihilates $N_0/R\text{Im}(\phi_{N_0}(N_0))$). As pointed out in [Va6, Rm. 6.4 (a)], the axioms [Va6, 6.2 (iii) to (vi)] hold automatically while the axioms [Va6, 6.2 (i) and (ii)] are implied by 6.1 (*) and 6.2 (**).  

{\bf Proof 2.} We follow [Va5, Lemma 2.5.2] to provide a second proof to (***) that does not rely on [Va6]. To prove that (***) holds  we can assume that we have $t_{z,\alpha}\in\scrT(N_0)$ and $v_{\alpha}\in\scrT(L^*\otimes_{\dbZ} \dbZ_{(p)})$ for all $\alpha\in\scrJ$. Let $Y_0$ be the affine $\Spec(R)$-scheme that parametrizes isomorphisms between $(N_0,(t_{\alpha})_{\alpha\in\scrJ})$ and $(L\otimes_{\dbZ} R,(v_{\alpha})_{\alpha\in\scrJ})$. The key point is to define (as everywhere in [Va1] to [Va8])
\medskip
 
{\bf (****)} $Y$ to be the {\it schematic closure} of $Y_{0,\eta}$ in $Y_0$!

\medskip

The group scheme $G_R:=G_{\dbZ_{(p)}}\times_{\dbZ_{(p)}} R$ acts on $Y$ from the left. From 6.2 (**) and [Ki4, proof of Prop. 1.3.4] one gets that for each  local ring $\scrR$ of $R$ which is a discrete valuation ring, $Y$ has valued points in a faithfully flat discrete valuating ring extension of $\scrR$. Using this and the fact that the schematic closures commute with pulls back via flat morphisms, one gets that $Y_U$ is a torsor under $(G_R)_U$. The rest is standard: based on [CTS] one gets that $Y_U$ is the trivial torsor and therefore $Y$ has $U$-valued points and thus also $\Spec(R)$-valued points. Therefore (***) holds and moreover one has:

\medskip
{\bf (*****)} there exists an isomorphism  $(N_0,(t_{1,\alpha})_{\alpha\in\scrJ})\arrowsim (L^*\otimes_{\dbZ} R,(v_{\alpha})_{\alpha\in\scrJ})$. 

\bigskip\noindent
{\bf 6.4. Step 3 (lift of the filtration).} The same as Subsection 4.3 (see [Ki4, Subsubsects. 1.1 and 1.5.8]). The only difference is that [Ki, Subsubsect. 1.5.7 and p. 19] got stuck to [Va1, Subsect. 5.4] and therefore it uses $\Rtil_e$ instead of $R_e$ for $p>2$ and it uses $R_e$ for $p=2$. 

\bigskip\noindent
{\bf 6.5. Steps 4 and 5 (deformation theory and parallel transport of Hodge cycles).} These steps are  the same as Subsection 4.5 (see [Ki4, Subsect. 1.5]), except that [Ki4, Subsect. 1.5] makes two modifications which do rely on Corollary 2.3.4 though this Corollary is not proved in [Ki4]. The modifications in [Ki4, Subsect. 1.5] are as follows: due to ``Claim 2.3.2'' (one can replace it by Corollary 2.3.4), loc. cit. does not require any appeal to Faltings result of Subsection 4.4 and for $p=2$ can use as well Zink deformation theory of displays. We do not believe in the ``parallel transport'' argument of the proof of [Ki4, Prop. 2.3.5] (at least it is way to sketchy). 

If $p\ge 3$, then [Ki4, Prop. 1.5.8] is a very particular case of [Va5, Thm. 5.2].${}^5$ $\vfootnote{5}{Loc. cit. is stated for $p\ge 3$ but it makes sense even for $p=2$ provided the $2$-rank of the $2$-divisible group over $V$ is $0$. In fact we have developed methods which even work in same cases when $p=2$ and the $2$-rank is arbitrary, see [Va7, Part I, Thm. 6.6].}$ 
\bigskip\noindent
{\bf 6.6. Comments.} The main merits of Method III are:

\medskip
{\bf (a)}  that it does not require a division into cases (to be compared with the end of Subsection 4.7); and 

{\bf (b)} one gets a variant of Milne conjecture for crystalline $\dbZ_p$-representations of the Galois group of $V[{1\over p}]$, provided in the \'etale $\dbZ_p$-context one has a reductive group scheme.

\medskip
 The main weak points of Method III are: {\bf (i)} it works only for $p>2$ and for very exceptional cases with $p=2$ and, in connection to Theorem 1.1, it can not get better or stronger results, {\bf (ii)} it is only a variation of the Methods I and II, {\bf (iii)} it can not prove any form of the Milne conjecture without the assumption that in the \'etale context one has a reductive group scheme over $\dbZ_p$ (in particular, it does not provide a new proof to the general results [Va4, Thm. 1.3 and Cor. 1.4] that pertain to $p$-divisible groups over Witt rings of perfect fields which are endowed with arbitrary families of crystalline tensors) {\bf (iv)}  (despite a better tool than Theorem 4.1) it has already accumulated numerous gaps, errors, and mistakes  (in [Ki3] it is also claimed that one has $\scrN=\scrN^{\text{cl}}$, Step 4 of [Ki3, Subsect. 4.2] is incorrect, see ``Claim 2.3.1'', see ``Claim 2.3.2'' and Remark 2.3.3 (d), see the false claim of Remark 2.3.3 (a), see Subsection 6.3, see Subsection 6.5, instead of quoting [Va2] it quotes the erroneous proof of [Mo, Lemma 3.6], etc.), {\bf (v)} it is mixed up (like one uses fully all the rings $V$, $R$, $R_e$, $\Rtil_e$, like one uses four deformation theories and two crystalline theories, etc.), and {\bf (vi)} it is purely written. 

\bigskip\noindent
{\bf 6.7. Remark.} Whatever part of [Ki4, Sect. 3] might be correct, in essence it is only a variation of the methods developed in [Va1, Subsects. 3.4, 6.1, and 6.2] and [Va3, Subsubsects. 2.4.2, 2.4.3, 4.3, and 5.1] for passing from integral canonical models of Shimura varieties of Hodge type to integral canonical models of Shimura varieties of abelian type.

%\newpage
\bigskip
\references{37}
{\nspace{

\Ref[Bo]
A. Borel,
\sl Linear algebraic groups,
\rm Grad. Texts in Math., Vol. {\bf 126}, Springer-Verlag, New York, 1991.

\Ref[BLR]
S. Bosch, W. L\"utkebohmert, and M. Raynaud,
\sl N\'eron models,
\rm Ergebnisse der Mathematik und ihrer Grenzgebiete (3), Vol. {\bf 21}, Springer-Verlag, Berlin, 1990.

\Ref[CTS]
J.-L. Colliot-Th\'el\`ene et J.-J. Sansuc,
\sl Fibr\'es quadratiques et composantes connexes r\'eelles,
\rm Math. Ann. {\bf 244} (1979), no. 2, pp. 105--134. 

\Ref[dJ]
J. de Jong,
\sl Homomorphisms of Barsotti-Tate groups and crystals in positive characteristic,
\rm Invent. Math. {bf 134}  (1998),  no. 2, pp. 301--333. Erratum:  Invent. Math.  {\bf 138}  (1999),  no. 1, p. 225. 

\Ref[De1]
P. Deligne,
\sl Travaux de Shimura,
\rm S\'eminaire  Bourbaki, 23\`eme ann\'ee (1970/71), Exp. No. 389, Lecture Notes in Math., Vol. {\bf 244}, pp. 123--165, Springer-Verlag, Berlin, 1971.

\Ref[De2]
P. Deligne,
\sl Vari\'et\'es de Shimura: interpr\'etation modulaire, et
techniques de construction de mod\`eles canoniques,
\rm Automorphic forms, representations and $L$-functions (Oregon State Univ., Corvallis, OR, 1977), Part 2,  pp. 247--289, Proc. Sympos. Pure Math., {\bf 33}, Amer. Math. Soc., Providence, RI, 1979.

\Ref[De3]
P. Deligne,
\sl Hodge cycles on abelian varieties,
\rm Hodge cycles, motives, and Shimura varieties, Lecture Notes in Math., Vol.  {\bf 900}, pp. 9--100, Springer-Verlag, Berlin-New York, 1982.

\Ref[DG]
M. Demazure, A. Grothendieck, et al., 
\sl Sch\'emas en groupes, Vol. {\bf II}, 
\rm Lecture Notes in Math., Vol. {\bf 152}, Springer-Verlag, Berlin-New York, 1970.

\Ref[Fa]
G. Faltings,
\sl Integral crystalline cohomology over very ramified
valuation rings,
\rm J. of Amer. Math. Soc. {\bf 12} (1999), no. 1, pp. 117--144.

\Ref[Fo]
J.-M. Fontaine,
\sl Représentations $p$-adiques des corps locaux. I. (French) [$p$-adic representations of local fields. I],
\rm  The Grothendieck Festschrift, Vol. II,  249--309, Progr. Math., 87, Birkhäuser Boston, Boston, MA, 1990.

\Ref[Ki1]
M. Kisin,
\sl Crystalline representations and $F$-crystals,
\rm Algebraic geometry and number theory,  pp. 459--496, Progr. Math., 253, Birkh\"auser Boston, Boston, MA, 2006.

\Ref[Ki2]
M. Kisin,
\sl  Modularity of 2-adic Barsotti-Tate representations,
\rm preprint

\Ref[Ki3]
M. Kisin,
\sl Integral canonical models of Shimura varieties -- notes for Journees Arithmetiques 2007
\rm 10 pages note, March 2007.

\Ref[Ki4]
M. Kisin,
\sl Integral models for Shimura varieties of abelian type,
\rm preprint, Nov. 5 (?), 2008 (45 pages dvi file shim.dvi).

\Ref[Ko]
R. E. Kottwitz,
\sl Points on some Shimura varieties over finite fields,
\rm J. of Am. Math. Soc. {\bf 5} (1992), no. 2, pp. 373--444.

\Ref[LR]
R. Langlands and M. Rapoport,
\sl Shimuravariet\"aten und Gerben, 
\rm J. reine angew. Math. 378 (1987), pp. 113--220.

\Ref[Mi]
J. S. Milne,
\sl The points on a Shimura variety modulo a prime of good
reduction,
\rm The Zeta functions of Picard modular surfaces, pp. 153--255, Univ. Montr\'eal, Montreal, Quebec, 1992.

\Ref[Mo]
B. Moonen,
\sl Models of Shimura varieties in mixed characteristics,
\rm Galois representations in arithmetic algebraic geometry (Durham, 1996),  267--350, London Math. Soc. Lecture Note Ser., 254, Cambridge Univ. Press, Cambridge, 1998.

\Ref[No] 
R. Noot, 
\sl Models of Shimura varieties in mixed characteristc, 
\rm J. Algebraic Geom. {\bf 5} (1996), no. 1, pp. 187--207.

\Ref[Va0]
A. Vasiu,
\sl Integral canonical models for Shimura varieties of Hodge type,
\rm Ph.D. thesis, Princeton University, 1994.

\Ref[Va1]
A. Vasiu,
\sl Integral canonical models for Shimura varieties of preabelian type,
\rm Asian J. Math. {\bf 3} (1999), no. 2, pp. 401--518.

\Ref[Va2]
A. Vasiu,
\sl A purity theorem for abelian schemes,
\rm Michigan Math. J. {\bf 52} (2004), no. 1, pp. 71--81.

\Ref[Va3]
A. Vasiu,
\sl Integral canonical models of unitary Shimura varieties,
\rm Asian J. Math. 12 (2008), no. 2, pp. 151--176.

\Ref[Va4]
A. Vasiu,
\sl Projective integral models of Shimura varieties of Hodge type with compact factors,
\rm J. Reine Angew. Math. {}\bf 618} (2008), pp. 51--75. 

\Ref[Va5]
A. Vasiu,
\sl A motivic conjecture of Milne,
\rm manuscript dated November 11, 2008,  54 pages available at http://arxiv.org/abs/math.NT/0308202.

\Ref[Va6]
A. Vasiu,
\sl Extension theorems for reductive group schemes,
\rm manuscript dated April 2006, available at 
http://arxiv.org/abs/math.NT/0406508.

\Ref[Va7]
A. Vasiu,
\sl Integral models in unramified mixed characteristic (0,2) of Shimura varieties of hermitian orthogonal Shimura varieties of PEL type, Parts I and II,
\rm http://arxiv.org/abs/math/0307205 and http://arxiv.org/abs/math/0606698.

\Ref[Va8]
A. Vasiu,
\sl Good reductions of Shimura varieties of hodge type in arbitrary unramified mixed characteristic, Parts I and II,
\rm math 0707.1668 and math 0712.1572.

\Ref[Zi]
T. Zink,
\sl Isogenieklassen von Punkten von Shimuramannigfaltigkeiten mit Werten in einem endlichen K\"orper,
\rm Math. Nachr. {\bf 112} (1983), pp. 103--124.

\enddocument